\documentclass{article}

\usepackage{amssymb, amsmath}
\usepackage{amsthm}
\usepackage{mathrsfs}
\usepackage{amscd}
\usepackage[active]{srcltx}
\usepackage{verbatim}
\usepackage{graphicx}

\usepackage[ruled]{algorithm2e}

\SetKwComment{Comment}{$\triangleright$\ }{}
\SetKwInOut{Input}{Input}
\SetKwInOut{Output}{Output}
 \usepackage{algpseudocode}

\usepackage[utf8]{inputenc}
\usepackage[T1]{fontenc}
\usepackage{lmodern}

\usepackage{hyperref} 
\usepackage{todonotes}
\usepackage{enumerate}
\usepackage{mathtools}
\usepackage{bbm}

\mathtoolsset{showonlyrefs,showmanualtags}

\newcommand\ud{\,\mathrm{d}}

\newcommand{\Ee}{\mathcal{E}}

\newcommand{\Oo}{\mathcal{O}}

\newcommand{\EE}{\mathbb{E}}

\newcommand{\NN}{\mathbb{N}}

\newcommand{\PP}{\mathbb{P}}

\newcommand{\RR}{\mathbb{R}}

\newcommand{\Var}{\mathrm{Var}}

\newcommand{\Lip}{\mathrm{Lip}}

\renewcommand{\epsilon}{\varepsilon}

\theoremstyle{plain}
\newtheorem{theorem}{Theorem}[section]
\theoremstyle{remark}
\newtheorem{remark}[theorem]{Remark}

\theoremstyle{plain}

\newtheorem{lemma}[theorem]{Lemma}
\newtheorem{proposition}[theorem]{Proposition}
\newtheorem{definition}[theorem]{Definition}

\numberwithin{equation}{section}

\allowdisplaybreaks

\title{Sequential construction of spatial networks with arbitrary degree sequence and edge length distribution}

\author{
Ivan Kryven\footnotemark[1]~
and
Rik Versendaal\footnotemark[1]~\textsuperscript{,}\footnotemark[2]
}

\footnotetext[1]{
Mathematical Institute, Utrecht University, P.O. Box 80010, 3508 TA Utrecht, The Netherlands, E-mail: \texttt{i.v.kryven@uu.nl}.
}
\footnotetext[2]{
Delft Institute of Applied Mathematics, Delft University of 
Technology, P.O. Box 5031, 2600 GA Delft, The Netherlands, E-mail: \texttt{r.versendaal@tudelft.nl}.
}

\date{\today}

\begin{document}

\maketitle

\begin{abstract}
Complex systems, ranging from soft materials to wireless communication, are often organised as random geometric networks in which nodes and edges evenly fill up the volume of some space. Studying such networks 
is difficult because they inherit their properties from the embedding space as well as from the constraints imposed on the network's structure by design, for example, the degree sequence. Here we consider geometric graphs with a given distribution for vertex degrees and edge lengths and propose a numerical method for unbiased sampling of such graphs. We show that the method reproduces the desired target distributions up to a small error asymptotically, and that is some boundary cases only a positive fraction of the network is guaranteed to possible to construct.\\

    
MSC2020 Classes: 05C80, 68W20, 60D05, 65D18
\end{abstract}

\tableofcontents

\newpage

\section{Introduction}
Models for complex networks are often constructed to incorporate geometric properties: each vertex is given a coordinate in some metric space and the probability for a pair of vertices to be connected depends on the metric distance between them. 
Complex networks for which the embedding space is explicitly given  are also called spatial networks 
\cite{barthelemy2018morphogenesis,penrose2003random}. For example, in the geometric random graph \cite{penrose2003random} the coordinates are sampled uniformly at random from a unit square and  the edges connect pairs of vertices closer apart than a given threshold. Geographical threshold graphs generalise the latter model by assigning vertices randomly distributed weights that affect the threshold \cite{bradonjic2008structure,bradonjic2007wireless}.
Such models tend to be useful in various modelling contexts and may even be informative  if there is no obvious metric space to consider. In the latter cases authors refer to the latent space which can be inferred from the network structure  \cite{boguna2021network}, or the feature space in the network completion problem \cite{kim2011network}. 

To mention examples, spatial networks are ubiquitous in soft materials and polymer gels, where vertices represent  correspondingly monomer units and granules \cite{torres2021effect,papadopoulos2018network}. Likewise, in wireless communication the vertices are receivers/transmitters distributed on a plane and the edges represent neighbours at a close distance \cite{weedage2021impact,bradonjic2007wireless}. Even if in both of the latter examples the number of neighbours is limited by the design factors, such as chemical valency  or number of radio channels per unit, it is the overall structure of the network that defines mechanical properties of the material and the throughput capacity of a wireless network. Such interests in the `bulk' properties of spatial networks motivate asymptotic analysis of large spatial networks having an a priori defined degree sequence. Besides, a similar paradigm occurs in spatial rumours spreading models \cite{janssen2017rumors} and  fluid mechanics in porous materials \cite{gu2022graphical}, and there is a potential for applications to spatial game theory \cite{nanda2017spatial} and spatial epidemiological models \cite{durrett1995spatial}, where specific degree distributions may need to be imposed. 

Constructing a general spatial network with a given degree distribution is a computationally difficult task. Even a simpler problem of constructing a non-spatial network that has a given degree sequence is already a known NP-hard problem with a long history of research \cite{greenhill2022generating}.  That being said, randomised algorithms achieve such construction in linear time for undirected \cite{BKS10} and directed  \cite{van2021fast} graphs when such networks are sufficiently sparse. One of the largest challenges in designing fast algorithms is to ensure that the produced graphs remain simple while they are sampled uniformly at random from all admissible possibilities that satisfy the constraints. The uniformity is an advantage when numerically studying  a property of the model as the ensemble average.
Furthermore, if the graphs are not required to be simple, a single iteration of stub matching will produce the desired multigraph and hence the construction problem is trivial in the non-spatial setting. In the spatial setting, however, it may happen that the algorithm fails to complete the whole graph because of incompatibility of the constraints, which is hard to assess a priori.

In this paper, we propose a numerical method that for a given set of points in space constructs a network with a target sequence for vertex degrees and distribution for the lengths of the edges. Controlling degree sequences independently of the spatial constraints opens up a possibility to single out the effect of space on the function of such networks and their bulk properties. In the rest of the paper we formalise our model for spatial networks and provide a sequential algorithm for constructing such networks numerically. The method we use to incorporate the target edge-length distribution may be useful for other sequential algorithms for graph sampling. 

In our analysis we show that the algorithm has asymptotically vanishing bias and study the convergence of the empirical distribution of the edge lengths to its target. It turns out that this convergence depends on how the target distribution scales with the number of vertices. In other words, controlling the density of the points in the metric space as a function of the number of vertices may affect the success rate of algorithm completion. We show that for some scaling regimes the convergence is achieved for arbitrary target distributions. At the same time, there is a boundary regime where the algorithm is guaranteed to complete only a linear fraction of edges in the graph. We provide a lower bound for this fraction when the edge length satisfies a certain dependency condition.


\section{The model and the sequential algorithm} \label{section:model}

\begin{figure}
    \centering
    \includegraphics[scale = 0.46]{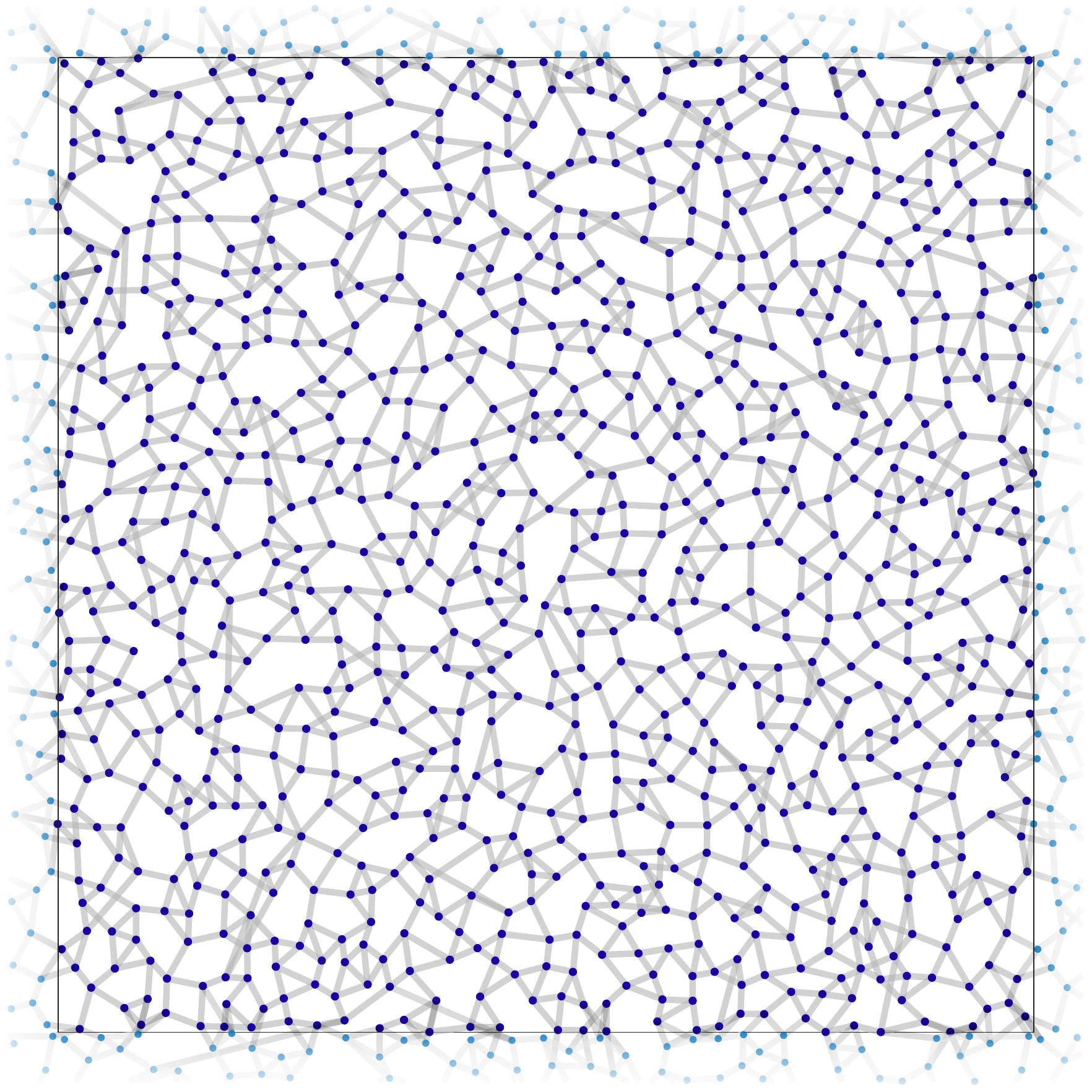}
    \caption{Spatial network with regular degree sequence and $1000$ vertices. The vertex coordinates are chosen on a torus  with Poisson disk model. The lengths of the edges are normally distributed with mean $0.03$ and variance $2 \cdot 10^{-5}$. The periodic boundaries are indicated with the box. }
    \label{fig:graph}
\end{figure}

\begin{algorithm}[h!]
\caption{Graph sampling with degrees and weights}
\hspace*{\algorithmicindent} \textbf{Input:} Graphical degree sequence $d$, target distribution $f$, and weights $\{r_{ij}\}_{1 \leq i < j \leq n}$. \\
    \hspace*{\algorithmicindent} \textbf{Output:} A graph $G$ with degree sequence $d^n$ and  edge-weight/length distribution, or 'failure'.
\begin{algorithmic}[1]
    \item Let $E$ be the set of edges;  $\hat d = (\hat d_1,\ldots, \hat d_n)$ a tuple of integers.\\ Initialize $E = \emptyset$ and $\hat d = d$.
    \item Choose two vertices $i, j\in V_n$ with probability proportional to 
    $$
    \frac{f_n(r_{ij})}{g_n(r_{ij})}\hat{d}_i\hat{d}_j\left(1-\frac{d_id_j}{4m}\right)
    $$
    among all pairs $i,j$ with $i\neq j$ and $ij \notin E$. Add $ij$ to $E$ and decrease $\hat{d}_i$ and $\hat{d}_j$ by 1.
    \item Repeat step 2 until no more edges can be added
    \item If $|E|< m$ report \emph{failure}, otherwise output 
    $G = (V,E)$.
\end{algorithmic}
\label{algorithm:sampling}
\end{algorithm}
We construct our spatial network on a set of vertices $V = \{1,\ldots,n\}$ by requiring three ingredients as input: 
\begin{enumerate}[1)]
    \item a graphical \cite{sierksma1991seven} degree sequence $d = (d_1,d_2,\dots,d_n)$,
    \item desired distribution of edge lengths, having density $f(x)$,
    \item a list of admissible edge lengths $r_{ij}\geq0$, with $r_{ij}=r_{ji}$ for $i,j=1,\dots,n.$ These numbers are identically distributed with density $g(x)$
and are pairwise independent.\\
\end{enumerate}

\begin{remark}[Spatial network]
 It is instructive to think of the vertices as  points randomly distributed in space, and of $r_{ij}$ as the distance between $i$ and $j$ ${i,j}$. Then the realisation of the graph satisfying the input constraints can be viewed as a spatial network, see, for example, Figure~\ref{fig:graph}.
Such spatial interpretation is only a little less general.
\end{remark}

Since $d_i$ is the degree of vertex $i$, the total number of edges in the network is $m = \frac12\sum_{i=1}^n d_i$.  Our aim is to randomly sample a simple graph $G=(V,E)$ satisfying the degree sequence $d$, such that the empirical distribution of edge lengths selected from $\{r_{ij}\}$ is close to the desired distribution $f(x)$. 
Algorithm~\ref{algorithm:sampling} starts with set $E$ being empty and sequentially adds edges one-by-one to $E$ while avoiding double edges or self loops. At each step, the next edge $\{i,j\}$ is selected with importance probability proportional to
$$
\frac{f_n(r_{ij})}{g_n(r_{ij})}\hat d_i\hat d_j\left(1 - \frac{d_id_j}{4m}\right).
$$
among all pairs $ij$ that maintain our graph being simple (\emph{i.e.} $i\neq j$  and $\{i,j\}\notin E$). Here, $\hat d_i$ is the remaining degree of vertex $i$ at the given step of the algorithm.  
Hence this probability is different at each step, which means that the probability of a given edge at step $r$ depends on the whole history of edges placed at earlier steps.

In the analysis section we show that if the degree sequences $d^n$ are uniformly bounded and $\frac{f_n}{g_n} = \Oo(n^\tau)$ for some $\tau < 1$, then Algorithm~\ref{algorithm:sampling} produces with high probability a graph with degree sequence $d^n$, and the empirical edge-weight distribution in this graph is close to the target distribution $f_n$.  Furthermore, we discuss an example showing that if $\tau = 1$, the result typically does not hold, and one has to terminate the algorithm early to assure the empirical edge-weight distribution of the graph is still close to the target distribution.

\subsection{Analysis of the algorithm} \label{section:main_result}

Let the random variables $E_k$ denote the edge chosen by Algorithm~\ref{algorithm:sampling} at iteration $k$ and random variables $L_k^n$ denote the weight of this edge. We then define the empirical weight distribution of the graph $G_n$, 
$$
L(G_n) = \frac1m\sum_{k=1}^m \delta_{L_k^n}.
$$
Here, $m = \frac12\sum_{i=1}^n d_i$ depends on $n$, and in particular,  $m \geq \frac12n$. Note that $L(G_n)$ is a random probability distribution. We will show that $L(G_n)$ converges to the target distribution $f$. For this we need to introduce an appropriate metric on the space of probability measures.

\begin{definition}\cite{Vil08}\label{def:Kantorovich}
Let $\mu,\nu$ be two probability measures on $\RR$. Let $||h||_{\Lip}$ denote the Lipschitz constant of a function $h:\RR \to \RR$. The Kantorovich distance $d_K(\mu,\nu)$ is defined as
$$
d_K(\mu,\nu): = \sup\left\{\int_\RR h \ud(\mu - \nu) \middle| h:M \to \RR \mbox{ continuous with } ||h||_{\Lip} \leq 1\right\}.
$$
\end{definition}


We are now ready to state the main theorem:
\begin{theorem}\label{theorem:convergence_edgelengths}
Let the degree sequences $d^n = (d_i)_{i \in V_n},$ be uniformly bounded and the weights $\{r_{ij}\}_{1 \leq i < j \leq n}$ be pairwise independent, identically distributed random variables with distribution $g_n$. Let $f_n$ be probability densities on $\RR_+$. Assume that $f_n,g_n$ all have support contained in some compact interval $I$. Assume that there exist constants $C_n > 0$ such that $f_n \leq C_ng_n$ for all $n \in \NN$ and $C_n = O(n^\tau)$ for some $\tau < 1$. Then for all $\epsilon > 0$, 
\begin{equation}\label{eq:convergence_edgelengths}
\lim_{n\to\infty} \PP\left(d_K(L(G_n),f_n) > \epsilon \right) = 0.
\tag{1}
\end{equation}
If $f_n$ converges to some $f$ in the Kantorovich distance, we can replace $f_n$ by $f$ in \eqref{eq:convergence_edgelengths}.
\end{theorem}



\subsection{Analysis of the boundary case $C_n = Cn$}

Our proof of Theorem~\ref{theorem:convergence_edgelengths} breaks down when $\tau = 1$. However, it turns out that this case is a very natural setting. To see this, replace the unit square with a square of arbitrary length, chosen such that the density of the points inside is constant for all $n$. Then $\tau = 1$ implies that the length of the edges in the graph is independent of $n$.  Another way to motivate this setting is to draw an analogy with the scaling typically used in random geometric graphs in the critical regime.

\subsubsection{Example: random geometric graphs}
\begin{figure}
    \centering
    \includegraphics[scale = 0.6]{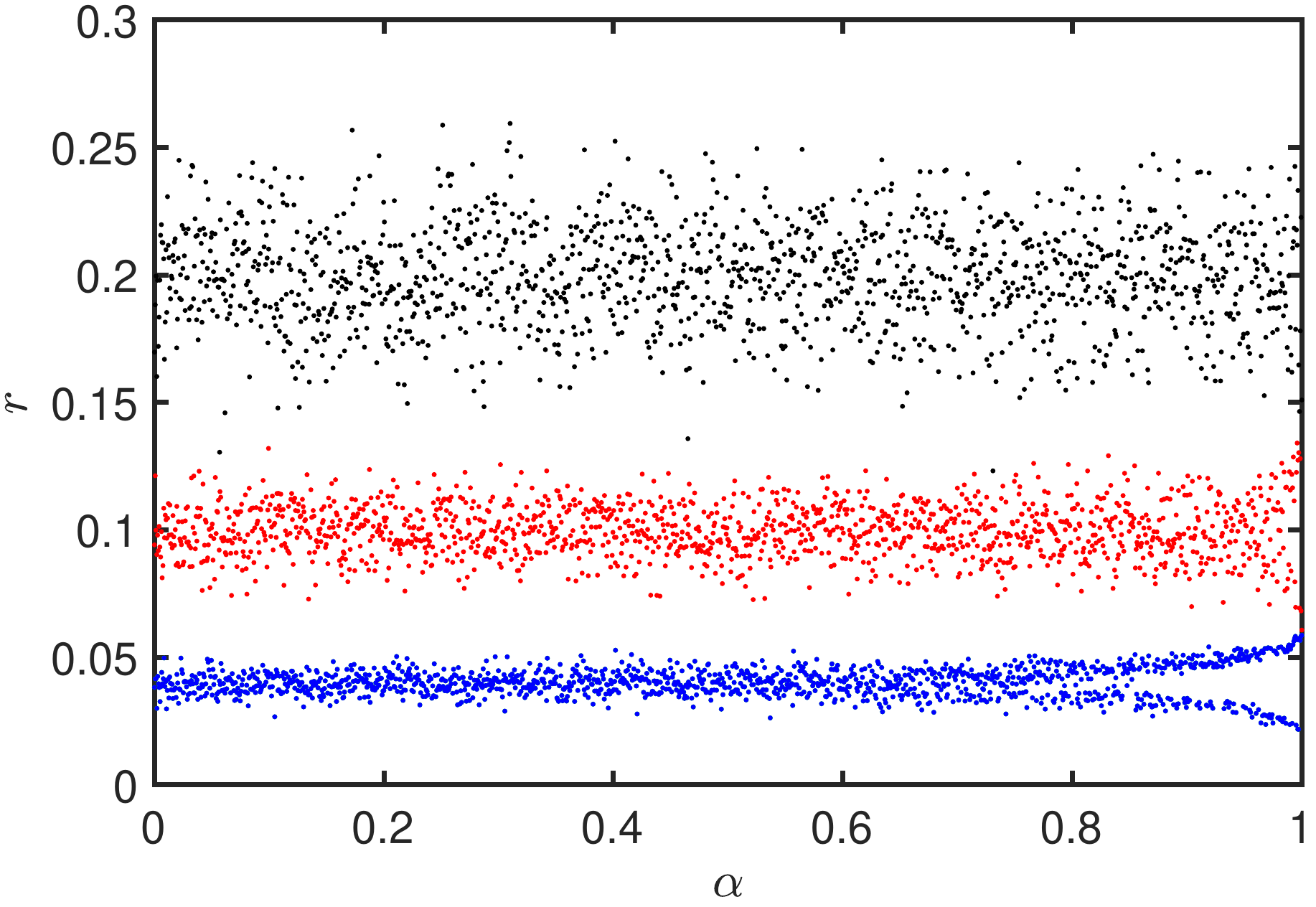}
    \caption{Length $r$ of edges placed by the algorithm at stage $\alpha$ is shown for three runs with different  edge length distributions.
    We use 1000 vertices with degree 3 and the target means being correspondingly $0.2$, $0.1$ and $0.04$ (top-down). The standard deviations are factor 0.15 times the means. The samples in the latter case split after $\alpha>0.8$ indicating that the target constraints are not compatible. }
    \label{figure:termination_fraction}
\end{figure}
Consider $n$ points $V_n = \{X_1,\ldots,X_n\}$ chosen independently and uniformly at random from the box $[0,1]^N$. We can construct the geometric graph $G(V_n,r)$ by connecting $X_i$ and $X_j$ whenever $|X_i - X_j| < r$. From classical percolation results (see e.g. \cite{Pen03}), it follows that the critical value of $r$ for the emergence of a giant component is of order $r_n \sim n^{-\frac1N}$.

Now suppose we want to sample such geometric graphs near this critical regime, under additional constraints on the degree distribution.  For this, let us set the weights $r_{ij} = |X_i - X_j|$, using periodic boundary conditions, i.e. we identify 0 with 1. By the triangle inequality, we find that any subset of these random variables is independent as long as the corresponding edges (in the complete graph) do not contain a cycle. More precisely, any two disjoint subsets of $\{r_{ij}\}$ are independent provided no cycle is formed containing edges from both subsets. In particular, the random variables are pairwise independent.
Now, to obtain a graph near the critical regime, our target distribution $f_n$ needs to be supported on an interval of the form $I_n = [0,\lambda n^{-\frac1N}]$. This implies that every vertex can on average connect only to $\lambda$ neighbours, \emph{i.e.}, the density of edges around each point is constant in $n$. Because the complete graph has $n-1$ edges per vertex, it follows that $\frac{f_n}{g_n} \sim n$, that is, $\tau = 1$.

Simulations shown in Figure~\ref{figure:termination_fraction}  suggest that in this case the edge-weight distribution of the graph produced by Algorithm~\ref{algorithm:sampling} fails to converge to the target distribution, while nevertheless remaining close to the target distribution for a positive fraction of the iterations. It turns out that this is true in general, provided we have a little more independence among the edge-weights than just pairwise. Let $G_n^k$ be the graph obtained by Algorithm~\ref{algorithm:sampling} after first $k$ iterations. We then have the following theorem.

\begin{theorem}\label{theorem:edge_weight_convergence_boundary_case}
Let the assumptions of Theorem~\ref{theorem:convergence_edgelengths} be satisfied, but instead assume that $C_n = Cn$ for some $C > 0$. Assume additionally that for any disjoint subsets together containing at most 4 edges, the corresponding sets of edge-weights are independent when no cycle is formed containing edges from both subsets. For $k \leq m$, define $G_n^k$ to be the random graph obtained by the algorithm after $k$ iterations. Then there is $\gamma > 0$ such that for every $\epsilon > 0$ we have
\begin{equation}\label{eq:convergence_edgelengths_boundary_case}
\lim_{n\to\infty} \PP(d_K(L(G_n^{\gamma m}),f_n) > \epsilon) = 0.
\tag{2}
\end{equation}
If $f_n$ converges to $f$ in the Kantorovich distance, we can replace $f_n$ by $f$ in \eqref{eq:convergence_edgelengths_boundary_case}.
\end{theorem}



\begin{remark}[Independence assumption in Theorem~\ref{theorem:edge_weight_convergence_boundary_case}]
The independence assumption in Theorem~\ref{theorem:edge_weight_convergence_boundary_case} implies that 3-tuples and 4-tuples of edge-weights may be dependent whenever they form a cycle. It is possible to  generalise these assumptions even further, as the proof requires that the amount of dependent 3-tuples and 4-tuples of edge-weights is sufficiently small.
\end{remark}



\section{Proof of Theorem~\ref{theorem:convergence_edgelengths}} \label{section:proof}

As before, we write $E_k$ for the (random) edge added by Algorithm~\ref{algorithm:sampling} at iteration $k$, and  $\Ee_k = \{E_1,\ldots,E_k\}$ for the set of edges formed at the first $k$ iterations. Furthermore, for $i = 1,\dots,n$ the remaining degree at iteration $k$ is denoted by
$$
d_i^k = d_i - \sum_{l=1}^{k-1} 1_{i \in E_l}.
$$
By definition of Algorithm~\ref{algorithm:sampling}, at iteration $k$, an edge $ij$ is chosen with probability
$$
p_{ij}^{k,n} = \frac{1}{Z_k^n}\frac{f_n(r_{ij})}{g_n(r_{ij})} d_i^kd_j^k\left(1 - \frac{d_id_j}{4m}\right)1_{ij \notin \Ee_{k-1}},
$$
where 
$$
Z_k^n 
= \sum_{ij} \frac{f_n(r_{ij})}{g_n(r_{ij})}d_i^kd_j^k\left(1 - \frac{d_id_j}{4m}\right)1_{ij \notin \Ee_{k-1}}
$$
is the normalization constant. Finally, for the sake of readability, in what follows we write 
$$
w_{ij} := 1 - \frac{d_id_j}{4m},
$$
since it is independent of step $k$.

\subsection{Sketch of the proof}
The proof relies on the idea that if the distribution of each $L_k^n$ is close to $f_n$, then so is their empirical distribution. Indeed, if $L_1^n,\ldots,L_m^n$ are i.i.d. random variables with distribution $f_n$, then standard results on laws of large numbers for measure-valued random variables (see e.g. \cite[Chapter 6]{DZ98}) show that $L(G_n)$ is close to $f_n$ in probability.
In our case, since the distribution of $L_k^n$ is given by conditioning on $\Ee_{k-1}$, the random variables are independent. However, they are not identically distributed. Therefore, we show in Proposition~\ref{prop:iteration_convergence} that with high probability the distribution of $L_k^n$ is close to $f_n$ in the Kantorovich distance. For this, we will make use of Chebyshev's inequality.
However, as seen in Proposition~\ref{prop:iteration_convergence}, this can only be achieved for all $k \leq m - o(m)$. This is because in order to obtain a vanishing bound from Chebyshev's inequality, the number of remaining edges at stage $k$ needs to tend to infinity when $n$ tends to infinity.
Once it is established that sufficiently many of the random variables $L_k^n$ have distribution close to $f_n$, the proof of Theorem~\ref{theorem:convergence_edgelengths}  follows from the technical result in Proposition~\ref{proposition:Wasserstein_LLN}, which is an adapted version of the law of large numbers for measure-valued random variables. 

\subsection{Convergence at each iteration}

As discussed in the previous section, we first show that at (almost) every iteration, the distribution of $L_k^n$ is close to $f_n$, uniformly in $k$. To prove this, we make use of Chebyshev's inequality. However, before we can do this, we first need a result showing that the sum of all weights remains large enough for sufficiently many iterations. For this, we first prove the following technical result. This is necessary to account for the fact that an edge not being chosen (yet) by the algorithm changes its distribution.

\begin{proposition}\label{prop:Tech_conditional_distribution_weights}
Let $c_1,\ldots,c_n > 0$ and let $X_1,\ldots,X_n$ be non-negative, pairwise independent random variables with distributions $g_1,\ldots,g_n$. Define the random variable $I \in \{1,\ldots,n\}$ with 
$$
\PP(I = i|X_1,\ldots,X_n) = \frac{c_iX_i}{\sum_{j=1}^n c_jX_j}.
$$
Then the distribution of $X_k$ conditional on $I \neq k$ is given by
$$
g_k(x|I\neq k) = \EE\left(\frac{\sum_{i\neq k}c_iX_i}{\sum_{j=1}^n c_jX_j}\right)^{-1}\EE\left(\frac{\sum_{i\neq k}c_iX_i}{c_kx + \sum_{i\neq k} c_iX_i}\right)g_k(x). 
$$
\end{proposition}
\begin{proof}
By Bayes' rule we have
$$
g_k(x|I\neq k) = \frac{\PP(I\neq k| X_k = x)g_k(x)}{\PP(I \neq k)}.
$$
Denoting by $g$ the joint distribution of $X_1,\ldots,X_n$, we have
\begin{align*}
\PP(I \neq k) 
&= 
\int \PP(I \neq k|X_1 = x_1,\ldots,X_n = x_n)g(x_1,\ldots,x_n) \ud(x_1,\ldots,x_n)
\\
&=
\int \frac{\sum_{i\neq k}c_ix_i}{\sum_{j=1}^n c_jx_j}g(x_1,\ldots,x_n) \ud(x_1,\ldots,x_n)
=
\EE\left(\frac{\sum_{i\neq k}c_iX_i}{\sum_{j=1}^n c_jX_j}\right).
\end{align*}
Since a similar expression holds for $\PP(I\neq k|X_k = x)$, the claim follows.
\end{proof}

An important consequence of Proposition~\ref{prop:Tech_conditional_distribution_weights} is that we can control the distribution of the weights $\frac{f_n(r_{ij})}{g_n(r_{ij})}$ conditional on the fact that edge $ij$ was not chosen by the algorithm yet. We need this to prove the following.

\begin{lemma}\label{lemma:sum_of_weights}
Let the assumptions of Theorem~\ref{theorem:convergence_edgelengths} be satisfied. Then there are constants\\ $C_1,C_2 > 0$ depending only on the maximum degree $d$, such that for all $\epsilon > 0$ and $\tau - \frac12 < \delta < \frac12$, 
$$
\PP(Z_k^n \geq (1-\epsilon)C_1n^{\frac32+\delta}|\Ee_{k-1}) \geq 1 -  \frac{C_n}{\epsilon^2C_2}n^{-\frac32-\delta}
$$
for all $k \leq m - m^{\frac34+\frac12\delta}$. Here, $C_n > 0$ and  $\tau < 1$ are such that $f_n \leq C_ng_n$  and $C_n = \Oo(n^\tau)$.
\end{lemma}
\begin{proof}
Denote by $R_k$ the set of valid edges, \emph{i.e.}
$$
R_k := \{ij|i \neq j, ij \notin \Ee_{k-1}, d_i^k,d_j^k > 0\}.
$$
We first find a lower bound for $|R_k|$. To this end, note that after $k \leq m - m^{\frac34+\frac12\delta}$ iterations, the total remaining degree is $2m - 2k \geq 2m^{\frac34 + \frac12\delta}$. Since the degree of every vertex is at most $d$, we have that at least $\frac{2m^{\frac34+\frac12\delta}}{d}$ vertices still have  left-over degrees. This gives rise to at least 
$$
\frac{1}{d^2}m^{\frac34+\frac12\delta}(m^{\frac34+\frac12\delta}-1) - (m - m^{\frac34+\frac12\delta})
$$
valid edges. Here, we subtract $k = m - m^{\frac34+\frac12\delta}$, since the algorithm has already chosen $k$ edges. Since furthermore $m \geq \frac12n$, we can find a constant $C_2 > 0$ such that
$
|R_k| \geq C_2n^{\frac32+\delta}.
$

Now consider $Z_k^n$ conditional on $\Ee_{k-1}$. Note that
$$
Z_k^n = \sum_{ij \in R_k} \frac{f_n(r_{ij})}{g_n(r_{ij})}d_i^kd_j^kw_{ij} \geq c\sum_{ij \in R_k} \frac{f_n(r_{ij})}{g_n(r_{ij})}
$$
for some $c > 0$. This is because for all $ij \in R_k$ we have $d_i^kd_j^k \geq 1$ and $w_{ij} \geq 1 - \frac{d^2}{4m}$.

As a consequence, we are done once we show that with sufficiently high probability
$$
\sum_{ij \in R_k} \frac{f_n(r_{ij})}{g_n(r_{ij})} \geq (1-\epsilon)C|R_k|.
$$
For this, let us write
$$
\mu_{ij}^k := \EE\left(\frac{f_n(r_{ij})}{g_n(r_{ij})}\middle| \Ee_{k-1}\right)
$$
so that
$$
\EE\left(\sum_{ij \in R_k} \frac{f_n(r_{ij})}{g_n(r_{ij})}\middle|\Ee_{k-1}\right) = \sum_{ij\in R_k} \mu_{ij}^k =: \mu_k.
$$

Using Chebyshev's inequality, we obtain
$$
\PP\left(\left|\sum_{ij \in R_k} \frac{f_n(r_{ij})}{g_n(r_{ij})} - \mu_k\right| > \epsilon\mu_k\middle|\Ee_{k-1}\right)
\leq \Var\left(\sum_{ij \in R_k} \frac{f_n(r_{ij})}{g_n(r_{ij})}\middle| \Ee_{k-1}\right)\epsilon^{-2}\mu_k^{-2}.
$$

Since the $\{r_{ij}\}$ are pairwise independent, also when conditioned on $\Ee_{k-1}$, we have
\begin{align*}
\Var\left(\sum_{ij \in R_k} \frac{f_n(r_{ij})}{g_n(r_{ij})}\middle|\Ee_{k-1}\right) 
&= 
\sum_{ij \in R_k} \Var\left(\frac{f_n(r_{ij})}{g_n(r_{ij})}\middle|\Ee_{k-1}\right) 
\\
&\leq 
C_n\sum_{ij \in R_k} \EE\left(\frac{f_n(r_{ij})}{g_n(r_{ij})}\middle| \Ee_{k-1}\right) 
\\
&= C_n\mu_k.
\end{align*}

Putting everything together, we find that
\begin{align*}
    \PP\left(\sum_{ij \in R_k} \frac{f_n(r_{ij})}{g_n(r_{ij})} \geq (1-\epsilon)\mu_k\middle|\Ee_{k-1}\right)
    &\geq 
    \PP\left(\left|\sum_{ij \in R_k} \frac{f_n(r_{ij})}{g_n(r_{ij})} - \mu_k\right| \leq \epsilon\mu_k\middle|\Ee_{k-1}\right)
    \\
    &\geq
    1 - \epsilon^{-2}C_n\mu_k^{-1}.
\end{align*}

The proof is completed by finding a suitable lower bound for $\mu_k$. We will show that
$$
\mu_k \geq \left(1 - \frac{2(k-1)C_n}{C_1\epsilon^2n^{\frac32+\delta}}\right)C_2n^{\frac32+\delta}
$$
for all $k \leq m - m^{\frac34+\frac12\delta}$ and $n$ being sufficiently large. We prove this using induction. 

For this, first note that $\mu_1 = \frac12n(n-1)$ is the total number of edges in $K_n$, which is certainly bigger than $C_2n^{\frac32+\delta}$ when $n$ is sufficiently large. Next, note that by Proposition~\ref{prop:Tech_conditional_distribution_weights} we have for $ij \in R_k$ that
\begin{align*}
\mu_{ij}^k
&\geq 
\EE\left(\frac{Z_{k-1}^n - \frac{f_n(r_{ij})}{g_n(r_{ij})}d_i^{k-1}d_j^{k-1}w_{ij}}{d^2C_n + Z_{k-1}^n - \frac{f_n(r_{ij})}{g_n(r_{ij})}d_i^{k-1}d_j^{k-1}w_{ij}}\middle| \Ee_{k-2}\right)\mu_{ij}^{k-1}
\\
&\geq 
\EE\left(\frac{Z_{k-1}^n - d^2C_n}{Z_{k-1}^n}\middle| \Ee_{k-2}\right)\mu_{ij}^{k-1}
\\
&\geq
\PP(Z_{k-1}^n \geq x|\Ee_{k-2})\left(1 - \frac{d^2C_n}{x}\right)\mu_{ij}^{k-1}.
\end{align*}

Applying this with $x = c(1-\epsilon)\mu_{k-1}$, we obtain
\begin{align*}
    \mu_k
    &\geq
    \left(1 - \frac{2(k-2)C_n}{C_2\epsilon^2n^{\frac32+\delta}}\right)\PP(Z_{k-1}^n \geq c(1-\epsilon)\mu_{k-1}|\Ee_{k-2})\left(1 - \frac{d^2C_n}{c(1-\epsilon)\mu_{k-1}}\right)|R_k|
    \\
    &\geq
    \left(1 - \frac{2(k-2)C_n}{C_2\epsilon^2n^{\frac32+\delta}}\right)\left(1 - \epsilon^{-2}C_n\mu_{k-1}^{-1}\right)\left(1 - \frac{d^2C_n}{c(1-\epsilon)\mu_{k-1}}\right)C_2n^{\frac32+\delta}
    \\
    &\geq
    \left(1 - \frac{2(k-2)C_n}{C_2\epsilon^2n^{\frac32+\delta}}\right)\left(1 - \left(\frac{d^2}{c(1-\epsilon)} + \frac{1}{\epsilon^2}\right)\frac{C_n}{C_2n^{\frac32+\delta}}\left(1 - \frac{2(k-2)C_n}{C_2\epsilon^2n^{\frac32+\delta}}\right)^{-1}\right)C_2n^{\frac32+\delta}
    \\
    &\geq
    \left(1 - \frac{2(k-2)C_n}{C_2\epsilon^2n^{\frac32+\delta}} - \frac{2C_n}{\epsilon^2C_2n^{\frac32+\delta}}\right)C_2n^{\frac32+\delta}
    =
    \left(1 - \frac{2(k-1)C_n}{C_2\epsilon^2n^{\frac32+\delta}}\right)C_2n^{\frac32+\delta}
\end{align*}
as desired.

Now, since $k \leq dn$, we find that
$$
\mu_k \geq  \left(1 - \frac{2C_n}{C_2\epsilon^2n^{\frac12+\delta}}\right)C_2n^{\frac32+\delta}
$$
for all $k$. Since furthermore $C_n = \Oo(n^\tau)$, we have for $\delta > \tau - \frac12$ that
$$
\lim_{n\to\infty} \left(1 - \frac{2C_n}{C_2\epsilon^2n^{\frac12+\delta}}\right) = 1.
$$
Therefore, we can find a constant $\tilde C_2 > 0$ such that for $n$ large enough we have
$$
\mu_k \geq \tilde C_2n^{\frac32+\delta}
$$
for all $k \leq m - m^{\frac34+\frac12\delta}$.

Collecting everything, we find
\begin{align*}
    \PP\left(Z_k^n \geq (1-\epsilon)c\tilde C_2n^{\frac32+\delta}\middle|\Ee_{k-1}\right)
    &\geq
    \PP\left(\sum_{ij \in R_k} \frac{f_n(r_{ij})}{g_n(r_{ij})} \geq (1-\epsilon)\mu_k\middle|\Ee_{k-1}\right)
    \\
    &\geq
    1 - \frac{C_n}{\epsilon^2\mu_k}
    \\
    &\geq
    1 - \frac{C_n}{\epsilon^2\tilde C_2}n^{-\frac32-\delta}
\end{align*}
as desired.
\end{proof}

Having a lower bound for $Z_k^n$, allows us to prove the following technical computation result.

\begin{lemma}\label{lemma:iteration_bound_start}
Let the assumptions of Theorem~\ref{theorem:convergence_edgelengths} be satisfied. Let $R_k$ be the set of valid edges at iteration $k$. Furthermore, define the quantities
$$
M_{ij}^k := \EE\left(\frac{f_n(r_{ij})}{g_n(r_{ij})}\middle| \Ee_{k-1}\right), \qquad \mbox{ and } \qquad  M_k = \sum_{ij \in R_k} M_{ij}^kd_i^kd_j^kw_{ij} 
$$
and
$$
H_{ij}^k := \EE\left(\frac{h(r_{ij})f_n(r_{ij})}{g_n(r_{ij})}\middle| \Ee_{k-1}\right) \qquad \mbox{ and } \qquad  H_k = \sum_{ij \in R_k} H_{ij}^kd_i^kd_j^kw_{ij}.
$$
Finally, suppose $0 \leq h(x) \leq b$ for all $x \in I$. Then for every $\epsilon > 0$ and all $\tau - \frac12 < \delta < \frac12$,
$$
\left|\frac{H_{ij}^k}{M_{ij}^k} - \int_I h(r)f_n(r)\ud r\right| < \epsilon
$$
for $n$ being sufficiently large and all $k \leq m - m^{\frac34 + \frac12\delta}$ and $ij \in R_k$. Moreover, we have
$$
\left|\frac{H_k}{M_k} - \int_I h(r)f_n(r) \ud r\right| < \epsilon
$$
for $n$ being sufficiently large.
\end{lemma}
\begin{proof}
We make use of Proposition~\ref{prop:Tech_conditional_distribution_weights}. Define the function
$$
F_k(r) = \EE_{k-1}\left(\frac{Z_{k-1,ij}^n}{d_i^{k-1}d_j^{k-1}w_{ij}r + Z_{k-1,ij}^n}\right),
$$
where
$$
Z_{k-1,ij}^n = Z_{k-1}^n - \frac{f_n(r_{ij})}{g_n(r_{ij})}d_i^{k-1}d_j^{k-1}w_{ij}.
$$
Note that 
$$
1 - F_k(r) = (d_i^{k-1}d_j^{k-1}w_{ij}r)\EE\left(\frac{1}{d_i^{k-1}d_j^{k-1}w_{ij}r + Z_{k-1,ij}^n}\right).
$$

Using Lemma~\ref{lemma:sum_of_weights}, we obtain
$$
1 - F_k(r) \leq \PP_{k-1}(Z_{k-1}^n \leq \eta n^{\frac32+\delta}) + \frac{1}{\delta n^{\frac32+\delta}} \leq \alpha n^{-\frac32-\delta}
$$
for some constant $\alpha > 0$ independent of $k$. 

Using this, together with Proposition~\ref{prop:Tech_conditional_distribution_weights} (see also the proof of Lemma~\ref{lemma:sum_of_weights}), we obtain
$$
\frac{H_{ij}^k}{M_{ij}^k} \leq \frac{(1 + \alpha n^{-\frac32-\delta})H_{ij}^{k-1}}{M_{ij}^{k-1}} = \left(1 + \frac{\alpha}{ n^{\frac32+\delta}}\right)\frac{H_{ij}^{k-1}}{M_{ij}^{k-1}}
$$
and
$$
\frac{H_{ij}^k}{M_{ij}^k} \geq \frac{H_{ij}^{k-1}}{(1 + \alpha n^{-\frac32-\delta})M_{ij}^{k-1}} = \left(1 - \frac{\alpha}{\alpha + n^{\frac32+\delta}}\right)\frac{H_{ij}^{k-1}}{M_{ij}^{k-1}}.
$$
Iterating this bound, and using that $k \leq m \leq dn$, we find that
$$
\left(1 - \frac{\alpha}{\alpha + n^{\frac32+\delta}}\right)^{dn}\frac{H_{ij}^1}{M_{ij}^1} \leq \frac{H_{ij}^k}{M_{ij}^k} \leq \left(1 + \frac{\alpha}{ n^{\frac32+\delta}}\right)^{dn}\frac{H_{ij}^1}{M_{ij}^1}. 
$$
The claim now follows by observing that the constants in the upper and lower bound converge to 1, and the fact that $M_{ij}^1 = 1$ and $H_{ij}^1 = \int_I h(r)f_n(r)\ud r$.

The second statement follows by rewriting the first to
$$
\left(\int_I h(r)f_n(r)\ud r - \epsilon\right)M_{ij}^k \leq H_{ij}^k \leq \left(\int_I h(r)f_n(r)\ud r + \epsilon\right)M_{ij}^k.
$$
and then summing over all $ij \in R_k$.
\end{proof}


We are now ready to prove that the distribution of $L_k^n$ (conditional on $\Ee_{k-1}$) is close to $f_n$.

\begin{proposition}\label{prop:iteration_convergence}
Let the assumptions of Theorem~\ref{theorem:convergence_edgelengths} be satisfied. Then for $n$ large enough there exists a constant $\alpha > 0$ only depending on $d$ and $|I|$ such that for every $\epsilon > 0$ small enough and $\delta < \frac12$ we have 
$$
\PP(d_K(L_k^n,f_n) > \epsilon|\Ee_{k-1}) \leq \frac{\alpha C_n}{\epsilon^2n^{\frac32 + \delta}}
$$
for all $k \leq m - m^{\frac34+\frac12\delta}$. Here, $C_n > 0$ are the constants such that $f_n \leq C_ng_n$.
\end{proposition}
\begin{proof}
Let $h$ be a Lipschitz function with $||h||_{\Lip} \leq 1$ and assume without loss of generality that $|I| \leq h(x) \leq 2|I|$ for all $x \in I$. It suffices to bound
$$
\PP\left(\left|\sum_{ij \in K_n} h(r_{ij})p_{ij}^{k,n} - \int_{\RR} h(r)f_n(r)\ud r\right| > \epsilon\middle|\Ee_{k-1}\right) 
$$
independent of $h$ and $k$.

For this, recall that given $\Ee_{k-1} = \{E_1,\ldots,E_{k-1}\}$ we have 
$$
p_{ij}^{k,n} = \frac{1}{Z_k^n}\frac{f_n(r_{ij})}{g_n(r_{ij})}d_i^kd_j^kw_{ij}1_{ij \notin \Ee_{k-1}}.
$$

Let $M_{ij}^k,M_k,H_{ij}^k,H_k$ and $R_k$ be as in Lemma~\ref{lemma:iteration_bound_start}. Then
$$
H_k = \EE\left(\sum_{ij \in R_k} \frac{h(r_{ij})f_n(r_{ij})}{g_n(r_{ij})}d_i^kd_j^kw_{ij}\middle|\Ee_{k-1}\right).
$$
Applying Chebyshev's inequality in a similar way as in the proof of Lemma~\ref{lemma:sum_of_weights} and using that $|I| \leq h(x) \leq 2|I|$ for all $x \in I$, we obtain
$$
    \PP\left(\left|\sum_{ij \in R_k} \frac{h(r_{ij})f_n(r_{ij})}{g_n(r_{ij})}d_i^kd_j^kw_{ij} - H_k \right| > H_k\epsilon\middle|\Ee_{k-1}\right)
    \leq
    \frac{C_n|I|d^2}{\epsilon^2H_k}.
$$

Likewise, we can prove 
$$
\PP\left(\left|Z_k^n - M_k\right| > M_k\epsilon\right)
\leq
\frac{C_nd^2}{2\epsilon^2M_k}.
$$

Now assume that
$$
\left|\sum_{ij \in R_k} \frac{h(r_{ij})f_n(r_{ij})}{g_n(r_{ij})}d_i^kd_j^kw_{ij} - H_k \right| \leq H_k\epsilon
$$
and
$
\left|Z_k^n - M_k\right| \leq M_k\epsilon.
$ 
Then we can estimate
\begin{align*}
    \MoveEqLeft
    \left|\sum_{ij \in R_k} h(r_{ij})p_{ij}^{k,n} - \int_{\RR} h(r)f_n(r)\ud r\right|
    \\
    &\leq
    \left|\frac{1}{Z_k^n}\sum_{ij \in R_k} \frac{h(r_{ij})f_n(r_{ij})}{g_n(r_{ij})}d_i^kd_j^kw_{ij} -\frac{H_k}{M_k}\right| + \left|\frac{H_k}{M_k} - \int_{\RR} h(r)f_n(r)\ud r\right|
    \\
    &\leq
    \frac{M_k}{Z_k^nM_k}\left|\sum_{ij \in R_k} \frac{h(r_{ij})f_n(r_{ij})}{g_n(r_{ij})}d_i^kd_j^kw_{ij} - H_k\right| + \frac{H_k|M_k - Z_k^n|}{Z_k^nM_k} + \epsilon
    \\
    &\leq
    \frac{2H_k\epsilon}{Z_k^n} + \epsilon.
\end{align*}

Here we used Lemma~\ref{lemma:iteration_bound_start} in the third line.

Now note that
\begin{align*}
H_k 
&= 
\sum_{ij \in R_k} \EE\left(\frac{h(r_{ij})f_n(r_{ij})}{g_n(r_{ij})}\middle|\Ee_{k-1}\right)d_i^kd_j^kw_{ij}
\\
&\leq
2|I|\EE\left(\frac{f_n(r_{ij})}{g_n(r_{ij})}\middle|\Ee_{k-1}\right)d_i^kd_j^kw_{ij}
\\
&= 
2|I| M_k.
\end{align*}

Furthermore, we have $Z_k^n \geq (1-\epsilon)M_k$. Using these estimates, we obtain
$$
\left|\sum_{ij \in R_k} h(r_{ij})p_{ij}^{k,n} - \int_{\RR} h(r)f_n(r)\ud r\right| \leq \frac{2|I|}{1-\epsilon}\epsilon + \epsilon. 
$$


Collecting everything, we find that there is a constant $\alpha_1 > 0$ such that for $\epsilon > 0$ small enough
$$
\left|\sum_{ij \in K_n} h(r_{ij})p_{ij}^{k,n} - \int_{\RR} h(r)f_n(r)\ud r\right| \leq \alpha_1\epsilon
$$
with probability at least
$$
1 - \frac{4C_n|I|d^2}{\epsilon^2H_k} - \frac{2C_nd^2}{\epsilon^2M_k}.
$$


Now note that since $h \geq |I|$, we have $H_k \geq |I|M_k$. Furthermore, we have
$$
M_k = \sum_{ij \in R_k} M_{ij}^kd_i^kd_j^kw_{ij} \geq c\sum_{ij \in R_k} M_{ij}^k = \mu_k \geq C_1n^{\frac32+\delta},
$$
where the latter follows from the proof of Lemma~\ref{lemma:sum_of_weights}. As a consequence, there exists a constant $\alpha >0$ such that the probability that
$$
\left|\sum_{ij \in K_n} h(r_{ij})p_{ij}^{k,n} - \int_{\RR} h(r)f_n(r)\ud r\right| > \epsilon
$$
is at most
$
\frac{\alpha C_n}{\epsilon^2n^{\frac32+\delta}}.
$
Because these estimates are independent of $h$, this completes the proof.
\end{proof}

\subsection{Converge of the edge-weight distribution}

Using Proposition~\ref{prop:iteration_convergence} we will show that the empirical edge-weight distribution $L(G_n)$ of $G_n$ is close to $f_n$. First, consider the following version of the law of large numbers for a sequence of random probability measures.

\begin{proposition}\label{proposition:Wasserstein_LLN}
Suppose $X_1^n,X_2^n,\ldots,X_n^n$ are real-valued independent random variables with random distributions $\mu_1^n,\ldots,\mu_n^n$. Assume that all supports of the $\mu_i^n$ are contained in a compact set $I$. Define
$$
\sigma_n = \frac1n\sum_{i=1}^n \delta_{X_i^n}.
$$
Furthermore, assume that there are distributions $\mu_n$ and for every $\epsilon > 0$ a sequence $\epsilon_n$ with $\lim_{n\to\infty} n\epsilon_n = 0$ such that
$$
\PP(d_K(\mu_i^n,\mu_n) > \epsilon) \leq \epsilon_n
$$
for all $i = 1,\ldots,n$ and all $n$ large enough. Then for all $\epsilon > 0$,
$$
\lim_{n\to\infty} \PP(d_K(\sigma_n,\mu_n) > \epsilon) = 0.
$$
\end{proposition}
\begin{proof}
Let $h$ be a Lipschitz function with $||h||_{\Lip} \leq 1$. Furthermore, let $X_n$ be distributed according to $\mu_n$. Then
$$
\int h\ud \sigma_n = \frac1n\sum_{i=1}^n h(X_i^n)
$$
and
$$
\int h\ud \mu_n = \EE(h(X_n)).
$$

Note that
\begin{align*}
    \PP\left(\left|\frac1n\sum_{i=1}^n h(X_i^n) - \EE(h(X_n))\right| > \epsilon\right)
    &\leq
    \PP\left(\left|\frac1n\sum_{i=1}^n h(X_i^n) - \frac1n\sum_{i=1}^n \EE(h(X_i^n))\right| > \frac\epsilon2\right)
    \\
    &\qquad +
    \PP\left(\left|\frac1n\sum_{i=1}^n \EE(h(X_i^n)) - \EE(h(X_n))\right| > \frac\epsilon2\right).
\end{align*}

We will bound both terms independent of $h$.
 Because the random variables are independent, the bound for the first term follows from Chebyshev's inequality,
$$
\PP\left(\left|\frac1n\sum_{i=1}^n h(X_i^n) - \frac1n\sum_{i=1}^n \EE(h(X_i^n))\right| > \frac\epsilon2\right) \leq \frac{4\sum_{i=1}^n \Var(h(X_i^n))}{n^2\epsilon^2}.
$$
Since $||h||_{\Lip} \leq 1$ and each random variable is contained in the compact set $I$, we may assume that $||h|_I||_\infty \leq |I|$. It follows that $\Var(h(X_i^n))$ can be uniformly bounded, independent of $h$. Hence, we find that
$$
\PP\left(\left|\frac1n\sum_{i=1}^n h(X_i^n) - \frac1n\sum_{i=1}^n \EE(h(X_i^n))\right| > \frac2\epsilon\right) \leq \frac{\tilde C}{n\epsilon},
$$
where $\tilde C$ is independent of $h$ and $n$.
For the second term, note that
\begin{align*}
    \PP\left(\left|\frac1n\sum_{i=1}^n \EE(h(X_i^n)) - \EE(h(X_n))\right| > \frac\epsilon2\right)
     &\leq
    \sum_{i=1}^n \PP\left(|\EE(h(X_i^n)) - \EE(h(X_n))| > \frac\epsilon2\right)
    \\
    &\leq
    \sum_{i=1}^n \PP\left(d_W(\mu_i^n,\mu_n) > \frac\epsilon2\right)
    \\
    &\leq
    n\epsilon_n.
\end{align*}

Collecting everything, we find that
$$
\PP\left(\left|\frac1n\sum_{i=1}^n h(X_i^n) - \EE(h(X_n))\right| > \epsilon\right) \leq \frac{\tilde C}{n\epsilon} + n\epsilon_n. 
$$
Since the upper bound is independent of $h$, and tends to 0, the claim follows.
\end{proof}

\begin{remark}\label{remark:Wasserstein_LLN}
Proposition~\ref{proposition:Wasserstein_LLN} remains true as long as the conditions are satisfied for all but at most $o(n)$ of the random variables. Since we divide by $n$, the influence of $o(n)$ random variables on the empirical average is negligible.
\end{remark}

\subsection{Proof of Theorem~\ref{theorem:convergence_edgelengths}}

With all the preparations done, the proof of Theorem~\ref{theorem:convergence_edgelengths} now follows from Proposition~\ref{prop:iteration_convergence} and Proposition~\ref{proposition:Wasserstein_LLN} (in particular Remark~\ref{remark:Wasserstein_LLN}).

\begin{proof}[Proof of Theorem~\ref{theorem:convergence_edgelengths}]
Given $\epsilon > 0$, it follows from Proposition~\ref{prop:iteration_convergence} that for $\tau -\frac12 < \delta < \frac12$ there exists a constant $\alpha > 0$ such that 
$$
\PP(d_K(L_k^n,f_n) > \epsilon) \leq \frac{\alpha C_n}{\epsilon^2n^{\frac32+\delta}}
$$
for all $k \leq m - m^{\frac34+\frac12\delta}$. Since $C_n = O(n^\tau) = o(n^{\delta + \frac12})$, we find that
$$
\lim_{n\to\infty} \frac{n\alpha C_n}{\epsilon^2n^{\frac32+\delta}} = \lim_{n\to\infty} \frac{\alpha C_n}{\epsilon^2 n^{\frac12 + \delta}} = 0.
$$
Hence, applying Proposition~\ref{proposition:Wasserstein_LLN}, especially Remark~\ref{remark:Wasserstein_LLN}, it follows that
$$
\lim_{n\to\infty} \PP(d_K(L(G_n),f_n) > \epsilon) = 0.
$$
This is the first claim. The second claim follows by an application of the triangle inequality.
\end{proof}



\section{Proof of Theorem~\ref{theorem:edge_weight_convergence_boundary_case}}\label{section:proof_boundary}

In this section, we provide a proof of Theorem~\ref{theorem:edge_weight_convergence_boundary_case}. One of the main difficulties compared to Theorem~\ref{theorem:convergence_edgelengths} is that in every iteration, the potential total weight of the edges that can no longer be considered by the algorithm can increase by $\Oo(n^2)$. In particular, this may result in the conditional expectation of the edge-weights to go to 0, making edges necessary to approximate the target edge-weight distribution unavailable.

The second difficulty is that due to the weights being $\Oo(n)$, Chebyshev's inequality no longer provides sharp enough concentration bounds. Therefore, using the additional independence assumptions, we will improve upon Chebyshev's inequality (see Lemma~\ref{lemma:Concentration_graph_dependence}).

\subsection{Preliminary estimates}

We start out by proving that the conditional expectation of the weights can be bounded both from above and away from 0. This is where conditions on the fraction $\gamma$ as in the statement of Theorem~\ref{theorem:edge_weight_convergence_boundary_case} arise. To do this, we first have to control the size of the total sum of the weights $Z_k^n$ using Chebyshev's inequality. Using these results, we then sharpen these estimates in the next section. The following is an adaptation  of Lemma~\ref{lemma:sum_of_weights}. 

\begin{lemma}\label{lemma:sum_of_weights_boundary_case}
Let the assumptions of Theorem~\ref{theorem:edge_weight_convergence_boundary_case} be satisfied. Then there exist constants $C_1,C_2 > 0$ (only depending on the maximum degree $d$) and a constant $\gamma^* < 1$ such that for all $\epsilon > 0$ and all $\gamma < \gamma^*$ we have
$$
\PP\left(Z_k^n \geq (1-\epsilon) C_1(1 - \gamma)^2n^2 \mbox{ for all } k \leq \gamma m\right) \geq 1 - \frac{C_2}{\epsilon^2n}.
$$
Here, $C > 0$ is such that $f_n \leq Cng_n$.
\end{lemma}
\begin{proof}
The proof is the same as that of Lemma~\ref{lemma:sum_of_weights}, apart from finding an appropriate lower bound for $\mu_k$. For this, suppose $k \leq \gamma m$ for $\gamma < \gamma^*$ to be determined later. Following the procedure of finding a lower bound for $\mu_k$ in Lemma~\ref{lemma:sum_of_weights}, we can prove that
$$
\mu_k \geq \left(1 - \frac{2(k-1)Cn}{L(\epsilon)|R_{\gamma m}|}\right)|R_{\gamma m}|,
$$
where 
$$
L(\epsilon) = \frac{d^2}{c(1-\epsilon)} + \frac{1}{\epsilon^2}.
$$

Furthermore, we can lower bound $|R_{\gamma m}|$. Indeed, after $\gamma m$ iterations, there is $2m - 2\gamma m = (1-\gamma)2m$ degree left. Since the maximum degree is bounded by $d$, this is distributed over at least $\frac{(1-\gamma)2m}{d}$ vertices. Therefore, we have 
$$
|R_{\gamma m}| \geq \eta\frac{(1-\gamma)^24m^2}{d^2}
$$
for some $\eta < 1$. 

Now write $\bar d = \frac1n\sum_{i=1}^n d_i$ for the average degree. Then $m = \frac12\bar dn$. Now, if $k \leq \gamma m$, we find
$$
\mu_k \geq \left(1 - \frac{\gamma Cd^2}{L(\epsilon)\bar d\eta(1-\gamma)^2}\right)\eta\frac{(1-\gamma)^2\bar d^2n^2}{d^2}
$$

The claim now follows if we take $\gamma^*$ such that
\begin{equation}
1 - \frac{\gamma^*Cd^2}{L(\epsilon)\bar d\eta (1-\gamma^*)^2} > 0.
\tag{3}
\label{eq:gamma_bound}
\end{equation}
\end{proof}

\begin{remark}
The condition for $\gamma^*$ in \eqref{eq:gamma_bound} only provides a lower bound for the threshold value. For large $n$, the constant $\eta$ will be close to 1. Furthermore, the function $L$ may be minimized over $\epsilon \in (0,1)$ to optimize the bound.
\end{remark}

We now show how to control the expectation of the weights at different iterations. This is necessary for sharpening the concentration inequality used in proving Lemma~\ref{lemma:sum_of_weights_boundary_case}

\begin{lemma}\label{lemma:bound_weigths}
Let the assumptions of Theorem~\ref{theorem:edge_weight_convergence_boundary_case} be satisfied. Define 
$$
M_{ij}^k := \EE\left(\frac{f_n(r_{ij})}{g_n(r_{ij})}\middle| \Ee_{k-1}\right).
$$
Then there exists a constant $\gamma^* \in (0,1]$ and constants $A,B > 0$ such that for all $k \leq \gamma^*m$ and all $ij \in R_k$ we have
$$
A \leq M_{ij}^k \leq B.
$$
\end{lemma}
\begin{proof}
We only consider the lower bound, the upper bound following similarly. Following the proof of Lemma~\ref{lemma:iteration_bound_start}, we have
$$
M_{ij}^k \geq \PP(Z_{k-1}^n \geq x|\Ee_{k-2})^{dn}\left(1 - \frac{d^2Cn}{x}\right)M_{ij}^{dn}.
$$
From Lemma~\ref{lemma:sum_of_weights_boundary_case}, by taking $x = \beta n^2$ for $\beta < 1$ small enough, we have
$$
M_{ij}^k \geq \left(1 - \frac{C_2}{n}\right)^{dn}\left(1 - \frac{d^2C\eta^{-1}}{n}\right)^{dn}
$$
Since the lower bound converges to $e^{-dC_2}e^{-d^3C\eta^{-1}} > 0$, the claim follows. 
\end{proof}

\subsection{A sharper concentration inequality}

In proving Lemma~\ref{lemma:sum_of_weights_boundary_case}, we only used that the edge-weights are pairwise independent. However, we assume slightly more independence. To strengthen the bound in Lemma~\ref{lemma:sum_of_weights_boundary_case}, we need the following improvement over Chebyshev's inequality. Because the proof of this sharper concentration inequality relies on Lemma~\ref{lemma:bound_weigths}, which in turn relies on Lemma~\ref{lemma:sum_of_weights_boundary_case}, we had to prove the weaker bound first.

\begin{lemma}\label{lemma:Concentration_graph_dependence}
Let $\{R_{ij}\}_{1\leq i < j \leq n}$ be random variables in $[0,Ln]$. Suppose the random variables are pairwise independent. Assume additionally that for any disjoint subsets together containing at most 4 edges, the corresponding sets of $R_{ij}$ are independent when no cycle is formed containing edges from different subsets. Define $\alpha_{ij} := \EE(R_{ij})$ and $\alpha = \sum_{ij} \alpha_{ij}$. Assume that there are constants $0 < c < C$ such that $c \leq \alpha_{ij} \leq C$ for all $ij$. Then there exists a constant $C_1 > 0$ such that for every $\epsilon > 0$ and $n \geq C$ we have
$$
\PP\left(\left|\sum_{ij} R_{ij} - \alpha\right| > \alpha\epsilon\right) \leq \frac{C_1}{\epsilon^4n^2}.
$$

\end{lemma}
\begin{proof}
By Markov's inequality we have
$$
\PP\left(\left|\sum_{ij} R_{ij} - \alpha\right| > \alpha\epsilon\right) \leq \frac{\EE\left(\left(\sum_{ij} R_{ij} - \alpha\right)^4\right)}{\epsilon^4\alpha^4}
$$

We now estimate the expectation. For this, first consider independent copies $\tilde R_{ij}$ which are all independent of each other. By expanding, we find that
\begin{align*}
\MoveEqLeft
\EE\left(\left(\sum_{ij} R_{ij} - \alpha\right)^4\right) - \EE\left(\left(\sum_{ij} \tilde R_{ij} - \alpha\right)^4\right)
\\
&=
\EE\left(\left(\sum_{ij} R_{ij}\right)^4\right) - \EE\left(\left(\sum_{ij} \tilde R_{ij}\right)^4\right)
\\
&\qquad + 4\alpha\EE\left(\left(\sum_{ij} R_{ij}\right)^3\right) - 4\alpha\EE\left(\left(\sum_{ij} \tilde R_{ij}\right)^3\right).
\end{align*}
Indeed, all lower order products cancel because of pairwise independence. This also implies that if we expand the powers, we only need to worry about terms with at least three different random variables. Next, define the sets
$$
C_3:= \{\{ij,jk,ki\}| 1 \leq i,j,k \leq n \mbox{ unequal}\}
$$
and
$$
C_4:= \{\{ij,jk,kl,li\}| 1 \leq i,j,k,l \leq n \mbox{ unequal}\}.
$$

Using the assumption on the 3-wise dependence of the random variables, we find that
\begin{align*}
\MoveEqLeft
\EE\left(\left(\sum_{ij} R_{ij}\right)^3\right) - \EE\left(\left(\sum_{ij} \tilde R_{ij}\right)^3\right)
\\
&= \sum_{(ij,jk,ki) \in C_3}  \EE(R_{ij}R_{jk}R_{ki}) - \EE(\tilde R_{ij})\EE(\tilde R_{jk})\EE(\tilde R_{ki})
\\
&= \sum_{(ij,jk,ki) \in C_3}  \EE(R_{ij}R_{jk}R_{ki}) - \EE(R_{ij})\EE(R_{jk})\EE(R_{ki}).
\end{align*}

Using that $R_{ij} \in [0,Ln]$, together with pairwise independence, we obtain 
$$
0 \leq \EE(R_{ij}R_{jk}R_{ki}) \leq Ln\EE(R_{jk})\EE(R_{ki})
$$
Using furthermore that $\EE(R_{jk}) \leq C$, it follows that for large enough $n$,
$$
\left|\sum_{(ij,jk,ki) \in C_3}  \EE(R_{ij}R_{jk}R_{ki}) - \EE(R_{ij})\EE(R_{jk})\EE(R_{ki})\right| \leq LC^2n|C_3| \leq LC^2n^4.
$$ Here we also used that $|C_3| \leq n^3$.

For the fourth power, this analysis is somewhat more involved. After expanding, and using the pairwise independence, we are left with
\begin{align*}
\MoveEqLeft
    \EE\left(\left(\sum_{ij} R_{ij}\right)^4\right) - \EE\left(\left(\sum_{ij} \tilde R_{ij}\right)^4\right)
    \\
    &=
    \sum_{(ij,jk,kl,li) \in C_4} \EE(R_{ij}R_{jk}R_{kl}R_{li}) - \EE(R_{ij})\EE(R_{jk})\EE(R_{kl})\EE(R_{li})
    \\
    &\qquad+
    \sum_{\substack{(ij,jk,ki) \in C_3\\ab \notin \{ij,jk,ki\}}}  \EE(R_{ab})(\EE(R_{ij}R_{jk}R_{ki}) - \EE(R_{ij})\EE(R_{jk})\EE(R_{ki}))
    \\
    &\qquad+
    \sum_{(ij,jk,ki) \in C_3}  \EE(R_{ij}^2R_{jk}R_{ki}) - \EE(R_{ij}^2)\EE(R_{jk})\EE(R_{ki}).
\end{align*}

Following a similar estimating procedure as before, we obtain
\begin{align*}
\left|\EE\left(\left(\sum_{ij} R_{ij}\right)^4\right) - \EE\left(\left(\sum_{ij} \tilde R_{ij}\right)^4\right)\right|
&\leq
L^2C^2n^2|C_4| + C^3Ln^3|C_3| + L^2C^2n^2|C_3|
\\
&\leq
(L^2C^2 + LC^3)n^6 + L^2C^2n^5.
\end{align*}
Here, we used in the last line that $|C_4| \leq n^4$.

Collecting everything, and noting that $\alpha \leq Cn^2$, we find that there exists a constant $C_1 > 0$, only depending on $C$ and $L$, such that
$$
\left|\EE\left(\left(\sum_{ij} R_{ij} - \alpha\right)^4\right) - \EE\left(\left(\sum_{ij} \tilde R_{ij} - \alpha\right)^4\right)\right| \leq C_1n^6.
$$

Using furthermore that $\alpha \geq \frac c2n(n-1)$, we thus obtain a constant $C_2 > 0$ such that
$$
\PP\left(\left|\sum_{ij} R_{ij} - \alpha\right| > \alpha\epsilon\right) \leq \frac{C_2}{\epsilon^4n^2} + \frac{\EE\left(\left(\sum_{ij} \tilde R_{ij} - \alpha\right)^4\right)}{\epsilon^4\alpha^4}
$$

Following a similar reasoning, we are done once we show that 
$$
\EE\left(\left(\sum_{ij} R_{ij} - \alpha\right)^4\right) \lesssim n^6.
$$

To this end, denote by $\kappa_l$ the $l$-th cumulant, i.e., the $l$-th coefficient of the power series of the log moment generating function. Then
\begin{align*}
\EE\left(\left(\sum_{ij} \tilde R_{ij} - \alpha\right)^4\right) 
&= \kappa_4\left(\sum_ij \tilde R_{ij}\right) + 3\Var\left(\sum_{ij} \tilde R_{ij}\right)^2
\\
&=
\sum_{ij}\kappa_4(\tilde R_{ij}) + 3\left(\sum_{ij} \Var(\tilde R_{ij})\right)^2
\\
&=
\sum_{ij} \EE((\tilde R_{ij} - \alpha_{ij})^4) + 3\sum_{ij\neq xy} \Var(\tilde R_{ij})\Var(\tilde R_{xy}).
\end{align*}
Here we used independence in the second line. Since $\tilde R_{ij} \in [0,Ln]$, we have $\EE((\tilde R_{ij} - \alpha_{ij})^4) \lesssim n^4$ so that
$$
\sum_{ij} \EE((\tilde R_{ij} - \alpha_{ij})^4) \lesssim n^6.
$$
Furthermore, we also have $\Var(\tilde R_{ij}) \lesssim n\EE(\tilde R_{ij}) \leq nC$, from which we obtain
$$
\sum_{ij\neq xy} \Var(\tilde R_{ij})\Var(\tilde R_{xy}) \lesssim n^6.
$$
This completes the proof.
\end{proof}

The following improvement of Lemma~\ref{lemma:sum_of_weights_boundary_case} is an immediate consequence of Lemma~\ref{lemma:Concentration_graph_dependence}.

\begin{lemma}\label{lemma:sum_of_weights_boundary_improved}
Let the assumptions of Theorem~\ref{theorem:edge_weight_convergence_boundary_case} be satisfied. There are constants $C_1,C_2 > 0$ (only depending on the maximum degree $d$) and a constant $\gamma^* < 1$ such that for all $\epsilon > 0$ and all $\gamma < \gamma^*$,
$$
\PP\left(Z_k^n \geq (1-\epsilon) C_1(1 - \gamma)^2n^2 \mbox{ for all } k \leq \gamma m\right) \geq 1 - \frac{C_2}{\epsilon^2n^2}.
$$
Here, $C > 0$ is such that $f_n \leq Cng_n$.
\end{lemma}
\begin{proof}
This is a direct consequence of replacing Chebyshev's inequality in the proof of\\ Lemma~\ref{lemma:sum_of_weights_boundary_case} by the inequality from Lemma~\ref{lemma:Concentration_graph_dependence}. This is allowed because of Lemma~\ref{lemma:bound_weigths}.
\end{proof}

The following result is the analogue of Lemma~\ref{lemma:iteration_bound_start}. The proof is identical, relying on the lower bound for $Z_k^n$ provided in Lemma~\ref{lemma:sum_of_weights_boundary_improved}.

\begin{lemma}\label{lemma:fraction_iteration_first}
Let the assumptions of Theorem~\ref{theorem:edge_weight_convergence_boundary_case} be satisfied. Let $R_k$ be the set of valid edges at iteration $k$. Furthermore, define the quantities $M_{ij}^k, M_{ij}, H_{ij}^k$ and $H_k$ as in Lemma~\ref{lemma:iteration_bound_start}. Finally, suppose $0 \leq h(x) \leq b$ for all $x \in I$. Then for every $\epsilon > 0$  and sufficiently large $n$,
$$
\left|\frac{H_{ij}^k}{M_{ij}^k} - \int_I h(r)f_n(r)\ud r\right| < \epsilon
$$
for all $k \leq \gamma^*m$ and $ij \in R_k$, where $\gamma^*$ is as in Lemma~\ref{lemma:sum_of_weights_boundary_improved}. Moreover, for all $k \leq \gamma^*$ and sufficiently large $n$ we have
$$
\left|\frac{H_k}{M_k} - \int_I h(r)f_n(r) \ud r\right| < \epsilon.
$$

\end{lemma}

\subsection{Proof of Theorem~\ref{theorem:edge_weight_convergence_boundary_case}}

With all preparations done, we are ready to prove Theorem~\ref{theorem:edge_weight_convergence_boundary_case}. For the sake of readability, we first prove that the edge-length distribution converges at each iteration separately. Theorem~\ref{theorem:edge_weight_convergence_boundary_case} is then an immediate consequence of Proposition~\ref{prop:iteration_convergence_boundary_case} together with Proposition~\ref{proposition:Wasserstein_LLN}.


\begin{proposition}\label{prop:iteration_convergence_boundary_case}
Let the assumptions of Theorem~\ref{theorem:edge_weight_convergence_boundary_case} be satisfied. Let $\gamma^*$ be as in Lemma~\ref{lemma:sum_of_weights_boundary_improved} and suppose $\gamma < \gamma^*$ and $k \leq \gamma m$. Then there exists a constant $C_1(C,d,|I|) > 0$ such that for every $\epsilon > 0$ small enough we have
$$
\PP(d_K(L_k^n,f_n) > \epsilon) \leq \frac{C_1(C,d,|I|)}{\epsilon^2n^2}.
$$
\end{proposition}
\begin{proof}
We follow the proof of Proposition~\ref{prop:iteration_convergence}. Assuming that $|I| \leq h(x) \leq 2|I|$, Lemma~\ref{lemma:bound_weigths} allows us to apply Lemma~\ref{lemma:Concentration_graph_dependence} instead of Chebyshev's inequality. From this, we obtain a constant $C_1 > 0$ such that with probability at least
$$
1 - \frac{C_1}{\epsilon^4n^2}
$$
we have
$$
\left|\sum_{ij \in R_k} \frac{h(r_{ij})f_n(r_{ij})}{g_n(r_{ij})} - H_k\right| \leq \epsilon H_k
$$
and
$$
|Z_k^n - M_k| \leq \epsilon M_k.
$$

Using these, we can estimate
\begin{align*}
\MoveEqLeft
    \left|\sum_{ij \in R_k} h(r_{ij})p_{ij}^{k,n} - \int_{\RR} h(r)f_n(r)\ud r\right| 
    \\
    &\leq
    \left|\frac1{Z_k^n}\sum_{ij \in R_k} \frac{h(r_{ij})f_n(r_{ij})}{g_n(r_{ij})}d_i^kd_j^kw_{ij} - \frac{H_k}{M_k}\right| + \left|\frac{H_k}{M_k} - \int_I h(r)f_n(r)\ud r\right|
    \\
    &\leq
    \frac{2\epsilon H_k}{Z_k^n} + \epsilon,
\end{align*}
where the last line holds for $n$ large enough by Lemma~\ref{lemma:fraction_iteration_first}. Now notice that
$$
Z_k^n \geq (1-\epsilon)M_k \geq \frac{1-\epsilon}{2|I|}H_k.
$$
Collecting everything, we have
$$
\left|\sum_{ij \in R_k} h(r_{ij})p_{ij}^{k,n} - \int_{\RR} h(r)f_n(r)\ud r\right| \leq \frac{4|I|}{1-\epsilon}\epsilon +  \epsilon
$$
with probability at least
$$
1 - \frac{C_1}{\epsilon^4n^2}.
$$
Since this is independent of $h$, the claim now follows. 
\end{proof}



\bibliographystyle{unsrt}
\bibliography{literature}

\end{document}